
\input amstex.tex
\documentstyle{amsppt}
\magnification1200
 \hsize=12.5cm
  \vsize=17.5cm
  \hoffset=1cm
  \voffset=2cm

\def\DJ{\leavevmode\setbox0=\hbox{D}\kern0pt\rlap
{\kern.04em\raise.188\ht0\hbox{-}}D}
\def\dj{\leavevmode
 \setbox0=\hbox{d}\kern0pt\rlap{\kern.215em\raise.46\ht0\hbox{-}}d}

\def\txt#1{{\textstyle{#1}}}
\baselineskip=13pt
\def\hf{{\textstyle{1\over2}}}
\def\b{\beta}
\def\e{\varepsilon}
\def\f{\varphi}
\def\G{\Gamma}
\def\k{\kappa}
\def\s{\sigma}
\def\t{\theta}
\def\={\;=\;}

\def\zt{\zeta(\hf+it)}

\def\D{\Delta}

\def\R{\Re{\roman e}\,} \def\I{\Im{\roman m}\,}
\def\z{\zeta}

\def\t{\theta}
\def\hf{{\textstyle{1\over2}}}
\def\txt#1{{\textstyle{#1}}}
\def\f{\varphi}

\def\le{\leqslant}
\def\ge{\geqslant}
\font\tenmsb=msbm10
\font\sevenmsb=msbm7
\font\fivemsb=msbm5
\newfam\msbfam
\textfont\msbfam=\tenmsb
\scriptfont\msbfam=\sevenmsb
\scriptscriptfont\msbfam=\fivemsb
\def\Bbb#1{{\fam\msbfam #1}}

\def \NN {\Bbb N}
\def \CC {\Bbb C}

\def \RR {\Bbb R}
\def \ZZ {\Bbb Z}
\def\d{\,{\roman d}}

\font\ff=cmr8
\def\txt#1{{\textstyle{#1}}}
\baselineskip=13pt

\font\teneufm=eufm10
\font\seveneufm=eufm7
\font\fiveeufm=eufm5
\newfam\eufmfam
\textfont\eufmfam=\teneufm
\scriptfont\eufmfam=\seveneufm
\scriptscriptfont\eufmfam=\fiveeufm
\def\mathfrak#1{{\fam\eufmfam\relax#1}}

\font\tenmsb=msbm10
\font\sevenmsb=msbm7
\font\fivemsb=msbm5
\newfam\msbfam
     \textfont\msbfam=\tenmsb
      \scriptfont\msbfam=\sevenmsb
      \scriptscriptfont\msbfam=\fivemsb
\def\Bbb#1{{\fam\msbfam #1}}

  \def\rightheadline{{\hfil{\ff
On the Rankin-Selberg problem in short intervals}
\hfil\tenrm\folio}}

  \def\leftheadline{{\tenrm\folio\hfil{\ff
  Aleksandar Ivi\'c }\hfil}}
  \def\emptyheadline{\hfil}
  \headline{\ifnum\pageno=1 \emptyheadline\else
  \ifodd\pageno \rightheadline \else \leftheadline\fi\fi}

\font\ff=cmr8
\font\teneufm=eufm10
\font\seveneufm=eufm7
\font\fiveeufm=eufm5
\newfam\eufmfam
\textfont\eufmfam=\teneufm
\scriptfont\eufmfam=\seveneufm
\scriptscriptfont\eufmfam=\fiveeufm
\def\mathfrak#1{{\fam\eufmfam\relax#1}}

\font\tenmsb=msbm10
\font\sevenmsb=msbm7
\font\fivemsb=msbm5
\newfam\msbfam
\textfont\msbfam=\tenmsb
\scriptfont\msbfam=\sevenmsb
\scriptscriptfont\msbfam=\fivemsb
\def\Bbb#1{{\fam\msbfam #1}}

\topmatter
\title
On the Rankin-Selberg problem in short intervals
\endtitle
\author
 Aleksandar Ivi\'c
\endauthor
\address
Katedra Matematike RGF-a, Universitet u Beogradu,  \DJ u\v sina 7,
11000 Beograd, Serbia.
\endaddress
\keywords the Rankin-Selberg problem, short intervals,
mean square  bounds
\endkeywords
\subjclass  11 N 37, 44 A 15, 26 A 12
\endsubjclass
\email
{\tt   ivic\@rgf.bg.ac.rs, aivic\_2000\@yahoo.com}
\endemail
\dedicatory
\enddedicatory
\abstract
If
$$
\D(x) \;:=\; \sum_{n\le x}c_n - Cx\qquad(C>0)
$$
denotes the error term in the classical Rankin-Selberg
problem, then we obtain a non-trivial upper bound for
the mean square of $\D(x+U) - \D(x)$ for a certain
range of $U = U(X)$. In particular, under the Lindel\"of hypothesis
for $\z(s)$, it is shown that
$$
\int_X^{2X} \Bigl(\D(x+U)-\D(x)\Bigr)^2\,{\roman d} x \;\ll_\e\; X^{9/7+\e}U^{8/7},
$$
while under the Lindel\"of hypothesis
for the Rankin-Selberg zeta-function the integral is bounded by $X^{1+\e}U^{4/3}$.
An analogous result for the discrete second moment of $\D(x+U)-\D(x)$ also holds.
\endabstract
\endtopmatter
\document
\head
1. Introduction and statement of results
\endhead
The classical Rankin-Selberg problem consists
of the estimation of the error term function
$$
\D(x) \;:=\; \sum_{n\le x}c_n - Cx,\leqno(1.1)
$$
where the notation is as follows. Let $\varphi(z)$ be a
holomorphic cusp form of weight $\kappa$ with respect to the full
modular group $SL(2,\ZZ)$, so that
$$
\f\left({az+b\over cz+d}\right) \;=\; (cz+d)^\k\f(z)\qquad\Bigl(a,b,c,d\in\ZZ,\;
 ad-bc=1\Bigr)
$$
when $\I z >0$ and $\lim_{\I z\to\infty}\f(z)=0$
(see e.g., R.A. Rankin [17] for basic notions).
We denote by $a(n)$ the $n$-th Fourier
coefficient of $\varphi(z)$ and suppose that $\varphi(z)$ is a
normalized eigenfunction for the Hecke operators $T(n)$, that is,
$a(1)=1$ and $T(n)\varphi=a(n)\varphi $ for every $n \in
\NN$. The classical example is $a(n) = \tau(n)$, when $\k=12)$.
This is the Ramanujan $\tau$-function defined by
$$
\sum_{n=1}^\infty \tau(n)x^n \=
x{\left\{(1-x)(1-x^2)(1-x^3)\cdots\right\}}^{24}\qquad(\,|x| < 1).
$$
The constant $C \,(\,>0)$ in (1.1) may be written down explicitly
(see e.g., [12]), and $c_n$ is the convolution function
defined by
$$
c_{n}\;=\;n^{1-\kappa}\sum_{m^2 \mid n}m^{2(\kappa-1)}
\left|a\Bigl({n\over m^2}\Bigr)\right|^2.\leqno(1.2)
$$
This is a multiplicative arithmetic function, namely $c_{mn} = c_mc_n$
when $(m,n)=1$, since $a(n)$ is multiplicative.
The classical Rankin-Selberg  bound of 1939 is
$$
\D(x) \;=\; O(x^{3/5}),\leqno(1.3)
$$
hitherto unimproved. In fact, this bound is one of the longest
standing unimproved bounds of Analytic number theory.
In their works, done independently, R.A.
Rankin [16] derives (1.3) from a general result of E. Landau [15],
while A. Selberg [19] states the result with no proof.
Note that, by the M\"obius inversion formula, (1.2) is equivalent
to
$$
|a(n)|^2n^{1-\k} \;=\; \sum_{d^2|n}\mu(d)c_{n/d^2}.
$$
Therefore using  (1.1), (1.3) and partial summation we obtain
$$
\sum_{n\le x}|a(n)|^2 \;= \; Dx^\k + O(x^{\k-2/5})\qquad(D>0),
$$
and conversely the above formula yields (1.1) with (1.3).

\medskip
Although
it seems very difficult at present to improve the bound in (1.3),
recently there have been some results on
the Rankin-Selberg problem (see the author's works [5]--[8]),
in particular on mean square estimates.
Namely, let as usual $\mu(\s)$ denote the Lindel\"of function
$$
\mu(\s) \;:= \;\limsup_{t\to\infty}\,{\log|\z(\s+it)|\over\log
t}\qquad\quad(\s\in\RR).\leqno(1.4)
$$
Then we have (see [6], [7]; the exponent of $\b$ was misprinted
as $2/(5-2\mu(\hf))$)
$$
\int_0^X \D^2(x)\d x \;\ll_\e\; X^{1+2\b+\e},\quad\b \= {2\over 5-4\mu(\hf)}.
\leqno(1.5)
$$
Here and later $\e$ denotes positive
constants which may be arbitrarily small, but are not necessarily the same at
each occurrence, while $\ll_\e$ means that the $\ll$--constant depends on $\e$.
Note that with the sharpest known  result (see M.N.
Huxley [2]) $\mu(\hf) \le 32/205$ we obtain $\b = 410/897 =
0.4457079\ldots\;$. The limit of (1.5) is the value $\b = 2/5$
if the Lindel\"of hypothesis for $\z(s)$ (that $\mu(\hf) =0$) is true.

\medskip
In this work we are interested in mean square bounds for $\D(x+U)-\D(x)
$ in the range $\;1\ll U\le X$,
especially when $U$ is ``short'', namely when $U = o(x)\;(x\to\infty)$.
First of all note that, since $c_n\ll_\e n^\e$, by (1.1) we  have
$$
\eqalign{
\D(x+U)-\D(x) &= \sum_{x<n\le x+U}c_n - CU\cr&
\ll_\e \sum_{x<n\le x+U} n^\e - CU \ll_\e Ux^\e.
\cr}\leqno(1.6)
$$
Although this bound may be considered as ``trivial", there does
not exist an analytic proof of it yet.
Hence using (1.5) and (1.6) we have
$$
\int_X^{2X}\Bigl(\D(x+U)-\D(x)\Bigr)^2\d x \ll_\e
\min\left(X^{1+2\b+\e}, X^{1+\e}U^2\right)
\quad(1 \ll U \le X).\leqno(1.7)
$$
One can call then (1.7) the ``trivial bound'' for the mean square of $\D(x+U)-\D(x)$,
and we seek a non-trivial bound, namely a bound which is (at least in certain
ranges of $U = U(X)$) sharper than (1.7).

\medskip
Recently there has been work on the analogue of this problem for some
related divisor problems. Let $\D_k(x)$ denote
the error term in the asymptotic formula for the summatory function of $d_k(n)$,
generated by $\z^k(s)\;(k\in\NN)$. Then in particular
$$
\D_2(x) = \sum_{n\le x}d(n) - x\bigl(\log x + 2\gamma -1\bigr)
\qquad\Bigl(d_2(n) \equiv d(n) = \sum_{\delta\mid n}1\Bigr)
$$
is the error term in the classical Dirichlet divisor problem
and $\gamma = -\G'(1)= 0.5772157\ldots\,$ is Euler's constant.
The author [9] proved that, for
$$
1 \ll U = U(X) \le \hf {\sqrt{X}}, \;c_3 = 8\pi^{-2}
$$
and computable constants $c_j$, we have
$$\eqalign{
\int_X^{2X}\Bigl(\D_2(x+U)-\D_2(x)\Bigr)^2\d x & = XU\sum_{j=0}^3c_j\log^j
\Bigl({\sqrt{X}\over U}\Bigr) \cr&
+ O_\e(X^{1/2+\e}U^2) + O_\e(X^{1+\e}U^{1/2}).\cr}\leqno(1.8)
$$
Thus for $X^\e \le U = U(X) \le X^{1/2-\e}$ it is seen that (1.8) is a true
asymptotic formula.

\medskip
A result analogous to (1.8) holds
if $\D_2(x+U)-\D_2(x)$ is replaced by the function
$E(x+U)-E(x)$, with different constants $c_j$, where
$$
E(T) := \int_0^T|\zt|^2\d t - T\Bigl(\log{T\over2\pi} + 2\gamma -1\Bigr)
$$
is the error term in the mean square formula for $|\zt|$.
For an extensive account on $E(T)$
see e.g., F.V. Atkinson's classical work [1], and the author's
monographs [3], [4].

\medskip
In the general case, when $k>2$, the above problem becomes more difficult. In [10] we
obtained mean square estimates for $\D_k(x+U) - \D_k(x)$. To formulate the results,
first we define $\s(k)$ as a number satisfying $\hf\le\s(k)<1$, for which
$$
\int_0^T|\z(\s(k)+it)|^{2k}\d t \;\ll_\e\;T^{1+\e}
$$
holds for a fixed integer $k\ge2$.
From zeta-function theory (see [3], and in particular Section 7.9  of
E.C. Titchmarsh [22]) it is
known that such a number exists for any given $k\in\NN$, but it is  not uniquely
defined, as one has
$$
\int_0^T|\z(\s+it)|^{2k}\d t \;\ll_\e\;T^{1+\e}\qquad(\s(k)\le\s<1).
$$
From Chapter 8 of [3] it follows that
 one has $\s(2) = \hf, \s(3) = {7\over12}, \s(4) = {5\over8},
\s(5) \le 9/20$ (see W. Zhang [23]) etc., but it is not
easy to write down (the best known
value of) $\s(k)$ explicitly as a function of $k$. Note that the  Lindel\"of
hypothesis that $\mu(\hf)=0$ is equivalent to the fact that
$\s(k) = \hf\; (\forall k \in \NN)$. Then the result
of [10] states: Let $k\ge3$ be a fixed integer. If $\s(k) = \hf$, then
$$
\int_X^{2X}{\Bigl(\D_k(x+U)-\D_k(x)\Bigr)}^2\d x \;\ll_\e\; X^{1+\e}U^{4/3}
\quad\Bigl(X^\e \le U = U(X) \le X^{1-\e}\Bigr).
\leqno(1.9)
$$
If $\hf < \s(k) < 1$, and $\t(k)$ is any constant satisfying $2\s(k)-1< \t(k)<1$,
then there exists $\e_1 = \e_1(k)>0$ such that
$$
\int_X^{2X}{\Bigl(\D_k(x+U)-\D_k(x)\Bigr)}^2\d x \ll_{\e_1} X^{1-\e_1}U^{2}
\quad\bigl(X^{\t(k)} \le U = U(X) \le X^{1-\e}\bigr).
$$
It is clear that, if the Lindel\"of hypothesis is true for $\z(s)$, then (1.9) holds
for all natural numbers $k \ge 2$.

\medskip
Recently the author and J. Wu [11] obtained a new upper bound for \break$\sum_{h\le H}\Delta_k(N,h)$
for $1\le H\le N$, $k\in \NN$, $k\ge3$,
where $\Delta_k(N,h)$ is the (expected) error term in the asymptotic formula
for $\sum_{N<n\le2N}d_k(n)d_k(n+h)$.

\bigskip
Now we state our results on the mean square of $\D(x+U)-\D(x)$ as the following

\medskip
THEOREM 1. {\it If $\mu = \mu(\hf)$ is defined by } (1.4) {\it then,
for $1\le U = U(X) \le X$,
$$
\int_X^{2X} \Bigl(\D(x+U)-\D(x)\Bigr)^2\d x \;\ll_\e\;
 X^{(9+12\mu)/(7+4\mu)+\e}U^{8/(7+4\mu)}.
\leqno(1.10)
$$
If
$$
Z(\hf +it) \;\ll_\e\;(|t|+1)^\e \leqno(1.11)
$$
holds, which is the Lindel\"of hypothesis for the Rankin--Selberg
zeta-function, then the above integral is bounded by $X^{1+\e}U^{4/3}$.}
\medskip
{\bf Corollary 1}. The bound in (1.10) improves (1.7) for
$$
X^{(1+4\mu)/(3+4\mu)} \;\le U \;\le X^{(16\mu^2-8\mu+9)/(20-16\mu)}.
$$

\medskip
{\bf Corollary 2}. If the Lindel\"of hypothesis for $\z(s)$ that
$\mu = \mu(\hf) = 0$ is true, then (1.10) reduces to
$$
\int_X^{2X} \Bigl(\D(x+U)-\D(x)\Bigr)^2\d x \;\ll_\e\; X^{9/7+\e}U^{8/7},
\leqno(1.12)
$$
and (1.12) improves (1.7) for $X^{1/3} \le U \le X^{9/20}$.

\medskip
There also exists a discrete analogue of Theorem 1. This is

\medskip
THEOREM 2. {\it If $\mu = \mu(\hf)$ is defined by } (1.4) {\it then,
for $1\le U = U(X) \le X$,
$$
\sum_{X<n\le{2X}} \Bigl(\D(n+U)-\D(n)\Bigr)^2 \;\ll_\e\;
 X^{(9+12\mu)/(7+4\mu)+\e}U^{8/(7+4\mu)}.
\leqno(1.13)
$$
If }(1.11) {\it holds, then the above sum is bounded by $X^{1+\e}U^{4/3}$.}
\medskip
It seems hard to ascertain what should be the true order of magnitude
of the function $\D(x+U)-\D(x)$. From (1.8) it seems plausible that
$$
\D_2(x+U) - \D_2(x) \ll_\e x^\e\sqrt{U}
\qquad\bigl(x^\e \le U = U(x) \le x^{1/2-\e}\bigr),
\leqno(1.14)
$$
which is a very strong conjecture made by M. Jutila [13], but it is not
clear whether there is sufficient  analogy between $\D(x+U)-\D(x)$
and $\D_2(x+U)-\D_2(x)$ to make any predictions about the order of
$\D(x+U)-\D(x)$ from (1.14).

\bigskip
\head 2. Proof of the Theorems
\endhead
There are two natural tools to study $\D(x)$. The first is the
explicit, truncated formula for $\D(x)$, of the Vorono{\"\i} type, namely
$$
\D(x) = {x^{3/8}\over2\pi}\sum_{k\le K}c_kk^{-5/8}\sin
\left(8\pi(kx)^{1/4}+{\txt{3\pi\over4}}\right) +
O_\e\Bigl(x^{3/4+\e}K^{-1/4}\Bigr),
$$
where the parameter $K$ satisfies $1 \ll K \ll x$. The proof of
this result can be found in [12]. However, the error term is much too large
for our present purpose. Therefore we resort to the use
of another natural tool in the study of $\D(x)$.
This is the Rankin--Selberg zeta-function
$$
Z(s) \;:=\; \sum_{n=1}^\infty c_n n^{-s},            \leqno(2.1)
$$
defined initially for $s = \s+it, \s > 1$, and for other values of
$s$ by analytic continuation. It has a simple pole at $s=1$ with residue
equal to $C$ (cf. (1.1)), and is otherwise regular.
For every $s\in\CC$ it satisfies the functional equation
$$
\G(s+\k-1)\G(s)Z(s) = (2\pi)^{4s-2}\G(\k-s)\G(1-s)Z(1-s).\leqno(2.2)
$$
The Rankin--Selberg zeta-function $Z(s)$ belongs to the
Selberg class ${\Cal S}$ of Dirichlet series of degree four.
For the definition and properties of
${\Cal S}$ see e.g., the seminal paper [20]
of A. Selberg and the review paper of Kaczorowski--Perelli [14].

\medskip
One also has the decomposition
$$
Z(s) := \sum_{n=1}^\infty c_n n^{-s} = \z(s)\sum_{n=1}^\infty b_n
n^{-s} = \z(s)B(s),\leqno(2.3)
$$
say, where $B(s)$ belongs to the  class ${\Cal S}$ of
of degree three, and moreover the function $B(s)$
is holomorphic for $\R s >0$. This follows from G. Shimura's work [21] (see
also A. Sankaranarayanan [18]). The coefficients $b_n$ in (2.3) are
multiplicative and satisfy $b_n \ll_\e n^\e$ (see [18]). Actually
the coefficients  $b_n$ are bounded by a log-power in mean square,
but this stronger property is not needed here.

\medskip
If we suppose that
$$
\int_X^{2X}|B( \hf+it)|^2\d t \;\ll_\e\; X^{\t+\e}\qquad(\t\ge
1),\leqno(2.4)
$$
and use the elementary fact (see Chapters 7 and 8 of [3] for the
results on the moments of $|\zt|$) that
$$
\int_X^{2X}|\zt|^2\d t \;\ll\; X\log X,\leqno(2.5)
$$
then from (2.3)--(2.5) and the Cauchy-Schwarz inequality for
integrals we obtain
$$
\int_X^{2X}|Z(\hf+it)|\d t \;\ll_\e\; X^{(\t+1)/2+\e}.\leqno(2.6)
$$
As $B(s)$ belongs to the Selberg class of degree three, then
$B(\hf+it)$ in (2.5) can be written as a sum of two Dirichlet
polynomials (e.g., by the reflection principle discussed in [3,
Chapter 4]), each of length $\ll X^{3/2}$, plus a manageable
error term. Thus by the mean value
theorem for Dirichlet polynomials (op. cit.) we have $\t \le 3/2$,
and any improvement on the value of $\t$ would give an improvement
of (1.3), as shown by the author in [6], [7].

\medskip

To prove (1.10), we start from (2.1) and Perron's inversion formula (see e.g., the
Appendix of [3]) to obtain
$$
\sum_{n\le x}c_n =
{1\over2\pi i}\int_{1+\e-i\tau}^{1+\e+i\tau}{x^s\over s}\,Z(s)\d s
+ O_\e(X^{1+\e}T^{-1}),\leqno(2.7)
$$
where $X\le x \le 2X,\, 1\ll \tau \ll X$ and $T \le\tau \le 2T$ will
be suitably chosen a little later.
We replace the segment of integration by the contour joining the points
$$
1+\e-i\tau,\, \hf-i\tau,\,\hf+i\tau,\,1+\e+i\tau.
$$
We encounter the simple pole of $Z(s)$ at $s=1$ of
and the residue will furnish $Cx$, the main term in (1.1).
Hence by the residue theorem (2.7) gives, once with $x$ and once with $x+U$,
$$
\eqalign{
\D(x+U) - \D(x)&= {1\over2\pi i}\int_{{1\over2}-i\tau}^
{{1\over2}+i\tau}{(x+U)^s-x^s\over s}\,
Z(s)\d s \cr&
\,+ O_\e(X^{1+\e}T^{-1}) + O(R(x,\tau)),\cr}
\leqno(2.8)
$$
where we set
$$
R(x,\tau) \;:=\; {1\over \tau}\int_{1\over2}^{1+\e}x^\s|Z(\s+i\tau)|\d \s.
$$
From (2.6) (with $\t = 3/2$) and the convexity of
mean values (see e.g., [3, Lemma 8.3]) we have
$$
\int_T^{2T}|Z(\s+it)|\d t \ll_\e T^{(3-\s)/2+\e}\qquad(\hf\le\s\le 1),
\leqno(2.9)
$$
and the integral in (2.9) is $\ll_\e T^{1+\e}$ for $\s\ge1$. It follows that
$$\eqalign{
\int_T^{2T}R(x,\tau)\d \tau &
\ll_\e {1\over T}\int_{1\over2}^{1+\e}x^\s\left(\int_T^{2T}|Z(\s+i\tau)|\d\tau\right)\d\s
\cr&
\ll_\e {1\over T}\max_{{1\over2}\le\s\le1+\e}\left({x\over\sqrt{T}}\right)^\s
T^{3/2+\e} \ll_\e XT^\e,\cr}
$$
since $T\ll X$. Note that this holds uniformly in $X \le x\le2X$.
Therefore there exists $T_0\in [T,\,2T]$ for
which
$$R(x,T_0) \;\ll_\e\; X^{1+\e}T^{-1}\qquad(X \le x\le2X)
$$
holds uniformly in $x$. It is $\tau = T_0$ that
is chosen in (2.7) and $T$ is the basic parameter to be determined.
Then using
$$
{(x+U)^s-x^s\over s} = \int_0^U (x+v)^{s-1}\d v
$$
we obtain from (2.8), since $T\le T_0\le2T$,
$$
\D(x+U) - \D(x)= {1\over2\pi i}\int_{{1\over2}-i\tau}^{{1\over2}+i\tau}
\left(\int_0^U (x+v)^{s-1}\d v\right)
Z(s)\d s + O_\e(X^{1+\e}T^{-1}).\leqno(2.10)
$$
On squaring (2.10) and integrating, we obtain
$$
\eqalign{&
\int_X^{2X}{\Bigl(\D(x+U)-\D(x)\Bigr)}^2\d x\cr&
\ll_\e \int_X^{2X}{\Bigl|\int_{-\tau}^\tau\int_0^U (x+v)^{{1\over2}-1+it}
Z(\hf+it)\d v\d t\Bigl|}^2\d x
+ X^{3+\e}T^{-2}.\cr}
\leqno(2.11)
$$
Let now $\psi(x)\;(\ge0)$ be a smooth function supported in $[X/2,\,5X/2]$, such that
$\psi(x) = 1$ when $X \le x \le 2X$ and $\psi^{(r)}(x) \ll_r X^{-r}\;(r=0,1,2,\ldots\,)$.
By using the Cauchy-Schwarz inequality for integrals it is seen that the integral
on the right-hand side of (2.11) does not exceed
$$
\eqalign{&
U\int_{X/2}^{5X/2}\psi(x)\int_0^U{\Bigl|\int_{-\tau}^\tau (x+v)^{-{1\over2}+it}
Z(\hf+it)\d t\Bigl|}^2\d v\,\d x\cr&
= U\int_0^U\int_{-\tau}^\tau\int_{-\tau}^\tau Z(\hf+it)Z(\hf-iy)J\,\d y\,\d t\,\d v,\cr}
$$
say, where
$$
J = J(X;v,t,y) :=\int_{X/2}^{5X/2}\psi(x)(x+v)^{-1}(x+v)^{i(t-y)}\d x.
$$
Integrating by parts we obtain, since $\psi(X/2) = \psi(5X/2) = 0$,
$$
J = {-1\over i(t-y)+1}\int_{X/2}^{5X/2}(x+v)^{i(t-y)}\Bigl(\psi'(x)
- {1\over x+v}\psi(x)\Bigr)\d x.
$$
By repeating this process it is seen that each time our integrand will be
decreased by the factor of order
$$
\ll \; {X\over|t-y|+1}\cdot {1\over X} \;\ll_\e\; X^{-\e}
$$
for $|t-y| \ge X^\e$. Thus if we fix any $A>0$, the contribution of
$|t-y|\ge X^\e$ will be $\ll X^{-A}$ if we integrate by parts $r = r(\e,A)$
times. For $|t-y|\le X^\e$ we estimate the corresponding contribution to $J$
trivially as $O(1)$ to obtain that the integral on the right-hand side of (2.11) is
$$
\eqalign{&
\ll_\e  U^2\int_{-\tau}^\tau\int_{-\tau,|t-y|\le X^\e}^\tau|Z(\hf+it)Z(\hf+iy)|
\d y\,\d t +1\cr&
\ll_\e  U^2\int_{-\tau}^\tau|Z(\hf+it)|^{2}\left(\int_{t-X^\e}^{t+X^\e}\d y\right)\d t
+1\cr&
\ll_\e U^2X^\e T^{{3\over2}+2\mu({1\over2})}.
\cr}
$$
Here we used  the elementary inequality
$
|ab| \le\hf\bigl(|a|^2 + |b|^2\bigr),
$
and the bound (cf. (2.4) with $\t = 3/2$)
$$
\eqalign{&
\int_X^{2X}|Z(\hf+it)|^2\d t = \int_X^{2X}|B(\hf+it)|^2|\zt|^2\d t \cr&
\ll_\e T^{2\mu({1\over2})+\e}\int_X^{2X}|B(\hf+it)|^2\d t
\ll_\e T^{{3\over2}+ 2\mu({1\over2})+\e}.\cr}
$$

Therefore it is seen that the left-hand side of (2.11) is
$$
\ll_\e\;X^\e(U^2T^{{3\over2}+2\mu({1\over2})}+ X^3T^{-2}).\leqno(2.12)
$$
With the choice
$$
T \;=\; X^{3/({7\over2}+2\mu({1\over2}))}U^{-2/({7\over2}+2\mu({1\over2}))}
$$
the terms in (2.12) are equalized. The condition
$1\ll T \ll X$ is trivial, and (2.12) yields (1.10).
Note that in proving (1.9) we could use power moments of $|\zt|$ for which there
is certainly more information than for the moments of $|Z(t)|$. This reflects
the quality of the bounds in (1.9) and (1.10).

\medskip

Finally note that if (1.11)
holds, which is the Lindel\"of hypothesis for the Rankin--Selberg
zeta-function, then obviously
$$
\int_X^{2X}|Z(\hf+it)|^2\d t \ll_\e X^{1+\e}.\leqno(2.13)
$$
This would replace (2.12) by
$$
\ll_\e X^\e(U^2T + X^3T^{-2}).
$$
The choice $T = XU^{-2/3}$  yields then
$$
\int_X^{2X}\Bigl(\D(x+U)-\D(x)\Bigr)^2\d x \;\ll_\e\; X^{1+\e}U^{4/3},\leqno(2.14)
$$
which is non-trivial in the whole range $1 \ll U \le X$.
Clearly for the proof (2.13) suffices instead of  the stronger (1.11).
The bound (2.14) is the
analogue of (1.9). This completes the proof of  Theorem 1.

\medskip
To prove (1.13) of Theorem 2 we employ the
method developed in [9]. We can assume that $U$
and $X$ are natural numbers, for otherwise we shall
make an admissible error by using trivial estimation.
Using (1.1) it is seen that integral  in (1.10) is equal to
$$
\eqalign{&
\sum_{X\le m\le 2X-1}\int_m^{m+1-0}\left(\sum_{x<n\le x+U}c_n -CU\right)^2\d x
\cr&
=  \sum_{X\le m\le 2X-1}\int_m^{m+1-0}\left(\sum_{m<n\le m+U}c_n -CU\right)^2\d x
\cr&
= \sum_{X\le m\le 2X-1}\int_m^{m+1-0}\left(\D(m+U)-\D(m)\right)^2\d x
\cr&
=\sum_{X\le m\le 2X-1}\left(\D(m+U)-\D(m)\right)^2\cr&
= \sum_{X\le n\le 2X}\left(\D(n+U)-\D(n)\right)^2 + O_\e(X^\e U).\cr}
$$
Here in the last step we used (1.6). Since the error term above is absorbed
in the expression on the right-hand side of (1.13), the proof of Theorem 2 is
finished.


\bigskip\bigskip
\Refs
\bigskip\bigskip

\item{[1]} F.V. Atkinson, The mean value of the Riemann zeta-function,
Acta Math. {\bf81}(1949), 353-376.

\item{[2]} M.N. Huxley, Exponential sums and the Riemann zeta-function V,
Proc. London Math. Soc. (4) {\bf90}(2005), 1-41.

\item{[3]} A. Ivi\'c,  The Riemann Zeta-Function,
   John Wiley \& Sons, New York, 1985 (2nd ed. Dover, Mineola, 2003).

\item{[4]} A. Ivi\'c, The mean values of the Riemann zeta-function,
LNs {\bf 82}, Tata Institute of Fundamental Research, Bombay (distr. by
Springer Verlag, Berlin etc.), 1991.

\item{[5]} A. Ivi\'c, Estimates of convolutions of certain
number-theoretic error terms, Intern. J. Math. and Math. Sciences,
Vol. 2004, No. {\bf1}(2004), 1-23.

\item{[6]} A. Ivi\'c, Convolutions and mean square estimates of certain
number-theoretic error terms, Publs. Inst. Math. {\bf80(94)}(2006), 141-156.

\item{[7]} A. Ivi\'c, On some mean square estimates in the
Rankin-Selberg problem, Applicable Analysis
and Discrete Mathematics {\bf1}(2007), 1-11.

\item{[8]} A. Ivi\'c, On the fourth moment in the Rankin-Selberg
problem,  Archiv Math. {\bf90}(2008), 412-419.

\item{[9]} A. Ivi\'c, On the divisor function and the Riemann zeta-function in
short intervals, The Ramanujan Journal, Volume {\bf19}, Issue 2 (2009), 207-224.

\item{[10]} A. Ivi\'c, On the mean square of the divisor function in short intervals,
Journal de Th\'eorie des Nombres de Bordeaux {\bf21}(2009), 195-205.

\item{[11]} A. Ivi\'c and J. Wu, On the general additive divisor problem,
see preprint at  arXiv:1106.4744.

\item{[12]} A. Ivi\'c, K. Matsumoto and Y. Tanigawa, On Riesz mean
of the coefficients of the Rankin--Selberg series, Math. Proc.
Camb. Phil. Soc. {\bf127}(1999), 117-131.

\item{[13]} M. Jutila, On the divisor problem for short intervals,
Ann. Univer. Turkuensis Ser. {\bf A}I {\bf186}(1984), 23-30.

\item{[14]} A. Kaczorowski and A. Perelli, The Selberg class: a
survey, in ``Number Theory in Progress, Proc. Conf. in honour
of A. Schinzel (K. Gy\"ory et al. eds)", de Gruyter,
Berlin, 1999, pp. 953-992.

\item{[15]}E.~Landau, \"{U}ber die Anzahl der
Gitterpunkte in gewissen Bereichen II, Nachr. Ges.
Wiss. G\"ottingen 1915, 209-243.

\item{[16]}R.~A.~Rankin, Contributions to the theory of Ramanujan's
   function $\tau(n)$ and similar arithmetical functions II. The order
   of the Fourier coefficients of integral modular forms,
    Proc. Cambridge Phil. Soc. {\bf 35}(1939), 357-372.

\item{[17]}R.~A.~Rankin, Modular forms and functions, Cambridge Univ. Press,
Cambridge, 1977.

\item{[18]} A. Sankaranarayanan, Fundamental properties of symmetric
square $L$-functions I, Illinois J. Math. {\bf46}(2002), 23-43.

    \item{[19]}A.~Selberg,  Bemerkungen \"{u}ber eine Dirichletsche Reihe,
    die mit der Theorie der Modulformen nahe verbunden ist,
    Arch. Math. Naturvid. {\bf 43}(1940), 47-50.

    \item{[20]} A.~Selberg, Old and new conjectures and results about a class
    of Dirichlet series, in ``Proc. Amalfi Conf. Analytic Number Theory 1989
    (E. Bombieri et al. eds.)",
    University of Salerno, Salerno, 1992, pp. 367--385.

    \item{[21]} G. Shimura, On the holomorphy of certain Dirichlet series,
    Proc. London Math. Soc. {\bf31}(1975), 79-98.

\item{[22]} E.C. Titchmarsh, The theory of the Riemann
zeta-function (2nd ed.),  University Press, Oxford, 1986.

\item {[23]} W. Zhang,
 On the divisor problem, Kexue Tongbao (in Chinese)
{\bf 33} (1988), 1484--1485.

\vskip1cm
\endRefs

\enddocument

\bye